\documentclass[11pt]{amsart}
\usepackage[OT4]{fontenc}
\usepackage{graphicx, amsmath, amssymb}
\newtheorem{Theorem}{Theorem}
\newtheorem{Corollary}[Theorem]{Corollary}

\newtheorem{Proposition}[Theorem]{Proposition}
\newtheorem*{example}{Example}
\theoremstyle{definition}
\newtheorem*{defn}{Definition}

\newcommand{\To}{\longrightarrow}
\newcommand{\so}{\Rightarrow}

\newcommand{\sU}{\mathcal{U}}
\newcommand{\sV}{\mathcal{V}}
\newcommand{\sW}{\mathcal{W}}
\newcommand{\sH}{\mathcal{H}}
\DeclareMathOperator{\as}{asdim} \DeclareMathOperator{\Prob}{Prob}
\DeclareMathOperator{\diam}{diam} \DeclareMathOperator{\mesh}{mesh}
\DeclareMathOperator{\Ind}{asInd} \DeclareMathOperator{\supp}{supp}
\DeclareMathOperator{\vol}{vol}
\begin{document}

\title{Asymptotic dimension in B\k{e}dlewo}%
\author{G. Bell}%
\address{Mathematical Sciences\\ UNC Greensboro \\
383 Bryan Building \\
Greensboro NC, 27402 USA}%
\email{gcbell@uncg.edu}%
\author{A. Dranishnikov}
\address{Department of Mathematics\\ University of Florida \\
PO Box 118105\\ 358 Little Hall \\ Gainesville, FL 32611 USA}
\email{dranish@math.ufl.edu}
\thanks{The second author was partially supported by NSF grant DMS-0305152}%


\bibliographystyle{amsplain}

\begin{abstract} This survey was compiled from lectures and problem
sessions at the International Conference on Geometric Topology at
the Mathematical Research and Conference Center in B\k{e}dlewo,
Poland in July 2005.
\end{abstract}
\maketitle

\section{Introduction}

These are the lecture notes from the Workshop in Asymptotic
Dimension Theory given at the International Conference on Geometric
Topology held in B\k{e}dlewo, Poland in July 2005.  The workshop
consisted of lectures on the basic theory by A. Dranishnikov and
problem sessions intended for young researchers hosted by G. Bell.

These notes seek to combine the lectures and problem sets into a
single survey which is intended to serve as a basic introduction to
the theory.

Section 2 of the notes gives the definition of asymptotic dimension
and useful equivalent formulations of the definition. A detailed
proof of the equivalences is given as well as some basic examples of
computations involving asymptotic dimension.  Also in this section
we prove that asymptotic dimension is a coarse invariant.

After proving a useful union theorem for asymptotic dimension, we
connect the asymptotic dimension to Lebesgue covering dimension via
the Higson corona.  Also in this section we introduce various other
notions of asymptotic dimension and give the known results on
coincidence of these notions.  The final part of section 2 is
devoted to uniform embedding of metric spaces with finite asdim into
products of trees.

The third section of the paper contains the main result of the
authors' paper \cite{BD3}: a Hurewicz-type theorem for asymptotic
dimension.  Although the proof of the theorem is omitted (it is
quite technical) the idea of the proof is given in the proof of a
much simpler version of the theorem that does not give a tight upper
bound.  The rest of this section is devoted to examples of
applications of the Hurewicz theorem to geometric group theory.

The other workshop at the Conference was given by M. Bestvina with
the assistance of L. Louder and H. Wilton and concerned limit
groups.  In the third section of the paper, we prove that limit
groups have finite asymptotic dimension.  This result follows from a
characterization of limit groups as constructible limit groups and
the Hurewicz-type theorem applied to the case of a group acting by
isometries on a tree.

The fourth section provides motivation for considering the
asymptotic dimension of a finitely generated group.  In particular,
we review results \cite{CG,DFW,Yu1,Yu2} concerning the Novikov
higher signature conjecture for finitely generated groups with
finite classifying spaces and finite asymptotic dimension.  In this
section we also give Higson and Roe's proof that finitely generated
groups with finite asymptotic dimension have G. Yu's property A,
i.e., are exact \cite{HR}.

The final section of the paper is devoted to the dimension function
of a finitely generated group.  In this section we give an example
of a finitely generated group with infinite asymptotic dimension and
show the quasi-isometry invariance of the growth of the dimension
function for a finitely generated group.

The authors wish to express their gratitude to the organizers of the
conference: Robert J. Daverman, Jerzy Dydak, Tadeusz Januszkiewicz,
Krystyna Kuperberg, S{\l}awomir Nowak, Stanis{\l}aw Spie\.{z}, and
chairman, Henryk Toru\'{n}czyk.  The conference was made possible
with the support of the Stefan Banach International Mathematical
Center, Warsaw University, Wroc{\l}aw University and the National
Science Foundation.

\section{Definitions and basic results}

\subsection{Definitions of asdim}

The asymptotic dimension can be defined for any coarse space (see
\cite{Ro03}), but we will consider only metric spaces in this paper.
We would like to view the asymptotic dimension of a space as somehow
dual to Lebesgue covering dimension. Given a cover $\sV$ of a
topological space, we say that the cover $\sU$ refines $\sV$ if
every $U\in\sU$ is contained in some element $V\in\sV.$  Recall that
the Lebesgue dimension $\dim$ of a topological space $X$ can be
defined as follows: $\dim X\le n$ if and only if for every open
cover $\sV$ of $X$ there is a cover $\sU$ of $X$ refining $\sV$ with
multiplicity $\le n+1.$  We will use the words order and
multiplicity of a cover interchangeably to mean the largest number
of elements of the cover meeting any point of the the space.  As
usual, we define $\dim X=n$ if it is true that $\dim X\le n$ but it
is not true that $\dim X\le n-1.$  We will do the analogous thing to
define $\as X=n.$

\begin{defn} Let $X$ be a metric space.  We say that the \emph{asymptotic
dimension} of $X$ does not exceed $n$ and write $\as X\le n$
provided for every uniformly bounded open cover $\sV$ of $X$ there
is a uniformly bounded open cover $\sU$ of $X$ of multiplicity $\le
n+1$ so that $\sV$ refines $\sU.$
\end{defn}

In practice one rarely uses this definition; instead one uses one of
the equivalent conditions described in the next theorem.  In
particular, for proving that the asymptotic dimension of a space is
finite it is easy to use condition (2) or (3).  In order to obtain a
tight upper bound for dimension, it is often better to work with
condition (5), which is phrased in terms of maps to uniform
complexes.

Before stating the theorem, we define the terminology used there.
Often we will need to consider very large positive constants and we
remind ourselves that they are large by writing $r<\infty$ instead
of $r>0.$  On the other hand, writing $\epsilon>0$ is supposed to
mean that $\epsilon$ is a small positive constant.

Let $r<\infty$ be given and let $X$ be a metric space. We will say
that a family $\sU$ of subsets of $X$ is $r$-{\em disjoint} if
$d(U,U')>r$ for every $U\neq U'$ in $\sU.$ Here, $d(U,U')$ is
defined to be $\inf\{d(x,x')\mid x\in U, x'\in U'\}.$  The $r$-{\em
multiplicity} of a family $\sU$ of subsets of $X$ is defined to be
the largest $n$ so that there is an $x\in X$ so that $B_r(x)$ meets
$n$ of the sets from $\sU.$  Recall that the {\em Lebesgue number}
of a cover $\sU$ of $X$ is the largest number $\lambda$ so that if
$A\subset X$ and $\diam(A)\le\lambda$ then there is some $U\in\sU$
so that $A\subset U.$

Let $K$ be a countable simplicial complex.  There are two natural
metrics one can place on $|K|$, the geometric realization of $K.$ We
wish to consider the uniform metric on $|K|.$  This is defined by
embedding $K$ into $\ell^2$ by mapping each vertex to an element of
an orthonormal basis for $\ell^2$ and giving it the metric it
inherits as a subspace.  A map $\varphi:X\to Y$ between metric
spaces is uniformly cobounded if for every $R>0$,
$\diam(\varphi^{-1}(B_R(y)))$ is uniformly bounded.

\begin{Theorem} Let $X$ be a metric space.  The following conditions
are equivalent.
\begin{enumerate}
    \item $\as X\le n$;
    \item for every $r<\infty$ there exist uniformly bounded,
    $r$-disjoint families $\sU^0,\ldots, \sU^n$ of subsets of $X$
    such that $\cup_i\sU^i$ is a cover of $X$;
    \item for every $d<\infty$ there exists a uniformly bounded
    cover $\sV$ of $X$ with $d$-multiplicity $\le n+1$;
    \item for every $\lambda<\infty$ there is a uniformly bounded
    cover $\sW$ of $X$ with Lebesgue number $>\lambda$ and multiplicity
    $\le n+1$; and
    \item for every $\epsilon>0$ there is a uniformly cobounded,
    $\epsilon$-Lipschitz map $\varphi:X\to K$ to a uniform simplicial
    complex of dimension $n.$
\end{enumerate}
\end{Theorem}

\begin{proof} This proof is not the most efficient one.  We prove (2)
$\so$ (3) $\so$ (4) $\so$ (5) $\so$ (2) and then (2) $\so$ (1) and
(1) $\so$ (4).

$(2)\so (3):$ Let $d<\infty$ be given and take $r>2d.$  We can find
uniformly bounded, $r$-disjoint families $\sU^0,\ldots,\sU^n$ of
subsets of $X$ covering $X.$  Put $\sV=\cup\sU^i.$ Suppose $x\in X$
and consider $B_d(x).$  If $U\cap B_d(x)$ and $U'\cap B_d(x)$ are
both nonempty, then $d(U,U')\le 2d<r,$ so $U$ and $U'$ must have
come from distinct families $\sU^i$ and $\sU^j$ with $i\neq j.$
Thus, there can be at most $n+1$ of the elements of $\sV$ having
non-trivial intersection with $B_d(x).$  So, (3) is proved.

$(3)\so(4):$  Let $\lambda<\infty$ be given and take a uniformly
bounded cover $\sV$ of $X$ with $5\lambda$-multiplicity $\le n+1.$
Define $\bar V=N_{2\lambda}(V).$ Then setting $\sW=\{\bar V\mid
V\in\sV\}$ we obtain a uniformly bounded cover with Lebesgue number
$>\lambda.$  It remains to show that the multiplicity is bounded by
$n+1.$  To this end, take $x\in X$ and observe that if $x\in \bar V$
then $d(x,V)< 2\lambda.$ So, $V$ is among the $n+1$ elements of
$\sV$ meeting $B_{2\lambda}(x).$  So the multiplicity of $\sW$ is
$\le n+1.$

$(4)\so (5):$ Let $\epsilon>0$ be given and suppose that $\sW$ is a
uniformly bounded cover of $X$ with multiplicity $\le n+1$ and
Lebesgue number greater than $\lambda=(2n+3)^2/\epsilon.$ For each
$W\in\sW,$ define the map $\varphi_W:X\to K$ by
\[\varphi_W(x)=\frac{d(x,X-W)}{\sum_{V\in\sW} d(x,X-V)}.\]
The maps $\{\varphi_W\}_W$ define a map $\varphi:X\to Nerve(\sW)$
where the nerve is a $n$-dimensional complex in $\ell^2$ with the
uniform metric.  It remains to check that $\varphi$ is uniformly
cobounded and $\epsilon$-Lipschitz.  If $\sigma$ is a simplex in $K$
and $x$ and $y$ both map to $\sigma,$ then there exist sets $U, V\in
\sW$ so that $x\in U$ and $y\in V$ and $U\cap V\neq \emptyset.$
Thus, $d(x,y)\le 2B,$ where $B$ is a uniform bound on the diameter
of the elements of $\sW.$

Finally, we check that $\varphi$ is $\epsilon$-Lipschitz. Let
$x,y\in X$ and $U\in\sW$. Let $\bar U$ denote the complement $X-U.$
The triangle inequality implies
\[|d(x,\bar U)-d(y,\bar U)|\le d(x,y).\] Also, observe that for any
$x\in X,$ $\sum_{U\in\sW}d(x,\bar U)\ge \lambda$ since $\lambda$ is
a Lebesgue number for $\sW.$ Thus we have
\[|\varphi_U(x)-\varphi_U(y)|=\left|\frac{d(x,\bar U)}
{\sum_{V\in\sW}d(x,\bar V)}-\frac{d(y,\bar
U)}{\sum_{V\in\sW}d(y,\bar V)}\right|\]
\[\le \frac{\left|d(x,\bar U)-d(y,\bar U)\right|}
{\sum_{V\in\sW}d(x,\bar V)}+\left|\frac{d(y,\bar U)}
{\sum_{V\in\sW}d(x,\bar V)}-\frac{d(y,\bar U)}
{\sum_{V\in\sW}d(y,\bar V)}\right|\]
\[\le\frac{d(x,y)}{\sum_{V\in\sW}d(x,\bar V)}+ \frac{d(y,\bar
U)}{\sum_{V\in\sW}d(x,\bar V)\sum_{V\in\sW}d(y,\bar
V)}\sum_{V\in\sW}\left|d(x,\bar V)- d(y,\bar V)\right|\]
\[\le\frac{1}{\lambda}d(x,y)+
\frac{1}{\lambda}\left(\sum_{V\in\sW}\left|d(x,\bar V)-d(y,\bar
V)\right|\right)\]
\[\le \frac{1}{\lambda}d(x,y)+\frac{2n+2}{\lambda}d(x,y)=\frac{(2n+3)}{\lambda}d(x,y).\]
Then we have
\[\|p(x)-p(y)\|_2=\left(\sum_{U\in\sW}|\varphi_U(x)-\varphi_U(y)|^2\right)^{\frac{1}{2}}\]\[\le
\left((2n+2)\left(\frac{(2n+3)}{\lambda}d(x,y)\right)^2\right)^{\frac{1}{2}}\]
\[\le\frac{(2n+3)^{3/2}}{\lambda}d(x,y)\le \epsilon d(x,y).\]

(5) $\so$ (2): Let $r$ be given and take $\phi:X\to K$ to be a
uniformly cobounded, $c/r$-Lipschitz map to a uniform complex of
dimension $n,$ where $c$ is a constant depending only on $n$ yet to
be determined. For each $i=0,\ldots,n,$ put
$\sV^i=\{St(b_\sigma,\beta^2K)\mid\sigma\subset K,\dim\sigma=i\},$
where $b_\sigma$ is the barycenter of $\sigma$ and $\beta^2K$
denotes the second barycentric subdivision.  Now, obviously there is
some $d$ so that $\diam V\le d$ for all $V\in\sV^i$ and there is a
constant $c$ depending only on the dimension of $K$ so that the
elements of $\sV^i$ are $c$-disjoint for all $i.$

Define $\sU^i=\{f^{-1}(V)\mid V\in\sV^i\}.$ Then $\diam(U)$ is
uniformly bounded since $f$ was uniformly cobounded.  Next since $f$
is $c/r$-Lipschitz, $d(U,U')<r$ implies that $d(f(U),f(U'))<c,$ so
the families $\sU^i$ are $r$-disjoint.

(2) $\so$ (1): Let $\sV$ be given with $\diam\sV\le \delta.$  Take
$r$-disjoint families $\sU^0,\ldots,\sU^n$ of uniformly bounded sets
in $X$ with $r>2\delta.$  Put $\bar\sU^i=\{N_\delta(U)\mid
U\in\sU^i\}.$  Put $\sU=\cup_i\bar\sU^i.$  Then since the
$\bar\sU^i$ are disjoint, the multiplicity of $\sU$ does not exceed
$n+1.$  Next, given $V\in\sV,$ $V$ must intersect some $U\in\sU^i$
for some $i.$  Since $\diam(V)\le\delta,$ $V\subset N_\delta(U)$
which is an element of $\sU.$  Thus, $\sV$ refines $\sU.$

(1) $\so$ (4): Let $\lambda<\infty$ be given.  Let
$\sV=\{B_\lambda(x)\mid x\in X\}.$  Clearly, $\sV$ is a cover of
$X$, so there is a uniformly bounded cover $\sU$ of $X$ with
multiplicity $\le n+1$ so that $\sV$ refines $\sU$.  Thus, any set
with diameter $\le\lambda$ will be entirely contained within one
element of $\sV$ so it will be entirely contained within one element
of $\sU,$ i.e., $L(\sU)\ge\lambda.$
\end{proof}

We conclude this section with a computation.

\begin{example} $\as T\le 1$ for all trees $T$ in the edge-length metric.
\end{example}

\begin{proof} Fix some vertex $x_0$ to be the root of the tree.  Let
$r<\infty$ be given and take concentric annuli centered at $x_0$ of
thickness $r$ as follows: $A_k=\{x\in T\mid
d(x,x_0)\in[kr,(k+1)r)\}.$  Although alternating the annuli (odd k,
even k) yields $r$-disjoint sets, these sets clearly do not have
uniformly bounded diameter.  We have to further subdivide each
annulus.

Fix $k>1.$ Define $x\sim y$ in $A_k$ if the geodesics $[x_0,x]$ and
$[x_0,y]$ in $T$ contain the same point $z$ with
$d(x_0,z)=r(k-\frac12).$ Clearly in a tree this forms an equivalence
relation.  The equivalence classes are $3r$ bounded and elements
from distinct classes are at least $r$ apart.  So, define $\sU$ to
be equivalence classes corresponding to even $k$ (along with $A_0$
itself) and $\sV$ to be equivalence classes corresponding to odd
$k.$  These two families cover $T$ and consist of uniformly bounded,
$r$-disjoint sets.  Thus, $\as T\le 1.$
\end{proof}

The proof that $\as T\le 1$ can be modified (see \cite{Ro05}) to
prove that $\delta$-hyperbolic metric spaces with bounded growth
have finite asymptotic dimension. Thus, finitely generated
$\delta$-hyperbolic groups have finite asdim, see also
\cite{BS,Gr93,Ro03}.

\subsection{Large-scale invariance of asdim}

One of the goals of this section is to prove that asymptotic
dimension is well-defined for finitely generated groups given the
word metric.

Let $\Gamma$ be a finitely generated group with finite, symmetric
generating set $S.$  We can define a norm on the group $\Gamma$
corresponding to $S$ by setting $\|\gamma\|_S$ equal to the minimal
number of $S$-letters necessary to present a word equal to $\gamma.$
Here we adopt the convention that the identity is presented by the
empty word. With this norm, we can define the (left-invariant) word
metric on $\Gamma$ by $d_S(g,h)=\|g^{-1}h\|_S.$  When $S$ is
understood we will simply write $d(g,h),$ see also Corollary
\ref{inv}

Let $X$ and $Y$ be metric spaces.  A map $f:X\to Y$ is a
$(\lambda,\epsilon)$-{\em quasi-isometry} if $d(f(x),f(x'))\le
\lambda d(x,x')+\epsilon$ for every pair of points $x,x'\in X$. A
map between metric spaces is a quasi-isometry if it is a
$(\lambda,\epsilon)$-quasi-isometry for some $\lambda>0$ and some
$\epsilon>0.$ The two spaces $X$ and $Y$ are quasi-isometric if
there is a quasi-isometry $f:X\to Y$ and a constant $C$ so that $Y
\subset N_C(f(X)).$

Coarse equivalence is a weaker notion of equivalence.  A map $f:X\to
Y$ between metric spaces is a {\em coarse embedding} if there exist
non-decreasing functions $\rho_1$ and $\rho_2$,
$\rho_i:\overline{\mathbb{R}}_+\to\overline{\mathbb{R}}_+$ such that
$\rho_i\to\infty$ and for every $x,x'\in X$
\[\rho_1(d_X(x,x'))\le d_Y(f(x),f(x'))\le
\rho_2(d_X(x,x'))\text{.}\]  Such a map is often called a {\em
coarsely uniform embedding} or just a {\em uniform embedding}.  The
metric spaces $X$ and $Y$ are {\em coarsely equivalent} if there is
a coarse embedding $f:X\to Y$ so that there is some $R$ such that
$Y\subset N_R(f(X)).$

Observe that quasi-isometric spaces are coarsely equivalent with
linear $\rho_i.$ (This is not entirely obvious from our definition,
one has to come up with a so-called quasi-inverse, which always
exists when spaces are quasi-isometric.) Also, one can always take
$\rho_2$ to be a linear function for coarse equivalence in finitely
generated groups.

\begin{Proposition} Let $f:X\to Y$ be a coarse equivalence.  Then
$\as X=\as Y.$
\end{Proposition}

\begin{proof} If $\sU^0,\ldots,\sU^n$ are $r$-disjoint, $D$-bounded
families covering $X$ then the families $f(\sU^i)$ are
$\rho_1(r)$-disjoint and $\rho_2(D)$-bounded.  Since $N_R(f(X))$
contains $Y$ we see that taking families $N_R(f(\sU^i))$ will cover
$Y$ and be $(2R+\rho_2(D))$-bounded and $(\rho_1(r)-2R)$-disjoint.
Since $\rho_i\to\infty,$ $r$ can be chosen large enough for
$\rho_1(r)-2R$ to be as large as one likes.  Therefore, $\as Y\le
\as X.$

The same proof applied to a coarse inverse for $f$ proves that $\as
X\le \as Y.$
\end{proof}

\begin{Corollary} \label{inv} Let $\Gamma$ be a finitely generated group.  Then
$\as\Gamma$ is an invariant of the choice of generating set, i.e.,
it is a group property.
\end{Corollary}

\begin{proof} Let $S$ and $S'$ be finite generating sets for
$\Gamma.$  We have to show that $(\Gamma,d_S)$ and $(\Gamma,d_{S'})$
are coarsely equivalent.  In fact, they are Lipschitz equivalent, as
we now show.

Let $\lambda_1=\max\{\|s\|_{S'}\mid s\in S\}$ and
$\lambda_2=\max\{\|s'\|_S\mid s'\in S'\}.$  It follows that
$\lambda_2^{-1}\|\gamma\|_{S'}\le\|\gamma\|_S\le
\lambda_2\|\gamma\|_{S'}.$  Take
$\lambda=\max\{\lambda_1,\lambda_2\}.$ Then $\lambda^{-1}
d_{S'}(g,h)\le d_S(g,h)\le \lambda d_{S'}(g,h).$
\end{proof}

\begin{example}  $\as\mathbb{R}=\as\mathbb{Z}=1.$
\end{example}

\begin{proof} First we show that $\mathbb{R}$ and $\mathbb{Z}$ are
coarsely isometric.  To this end, let $f:\mathbb{Z}\to\mathbb{R}$ be
the identity map.  Since the metric on $\mathbb{R}$ restricted to
$\mathbb{Z}$ is the same as the metric on the image of $f,$ this map
is a coarse equivalence.  Observe too, that
$N_1(f(\mathbb{Z}))=\mathbb{R}.$

We'll show that $\as\mathbb{Z}\le 1$ and $\as\mathbb{Z}\ge 1.$
First, given an $R<\infty$, consider the sets $A_k=[2Rk,2R(k+1)],$
where $k\in\mathbb{Z}.$  Let $\sU=\{A_{2k}\mid k\in \mathbb{Z}\}$
and $\sV=\{A_{2k-1}\mid k\in\mathbb{Z}\}.$  Clearly, the elements of
each of these families have diameter bounded by $2R.$  It is also
easy to check that any two elements from the same family are
$R$-disjoint as required.

On the other hand, if $\as\mathbb{Z}\le 0$ then for any $R<\infty$
there would be a cover of $\mathbb{Z}$ by an $R$-disjoint family of
subsets of $\mathbb{Z}$ with uniformly bounded diameter.  Let $R>1$
and let $U$ be the element of the family containing $0.$  Since
$d(n,n+1)=1$ no consecutive integers can belong to different
elements of $\sU.$  Thus, $\mathbb{Z}\subset U,$ so
$\diam(U)=\infty.$  Thus, $\as\mathbb{Z}>0.$
\end{proof}

The technique used in the proof that $\as\mathbb{Z}>0$ can also be
used to show the following result.

\begin{example} Let $\Gamma$ be a finitely generated group.  Then
$\as\Gamma=0$ if and only if $\Gamma$ is finite.
\end{example}

J. Smith \cite{Smi} has classified all countable (not just finitely
generated) groups with a left-invariant, proper metric and
asymptotic dimension $0$ and so is able to deduce this last example
as a corollary.

\begin{Theorem} Let $G$ be a countable group.  Then $\as G=0$ if and
only if every finitely generated subgroup of $G$ is finite.
\end{Theorem}

\begin{Proposition} Let $X$ be a metric space and $Y\subset X.$
Then $\as Y\le \as X.$
\end{Proposition}

\begin{proof} Let $R<\infty$ be given and take a cover $\sU$ of $X$
by uniformly bounded sets with $R$-multiplicity $\le n+1.$  Clearly
the restriction of this cover to $Y$ yields a cover whose elements
have uniformly bounded diameter and at most $n+1$ of them can meet
any ball of radius $R$ in $Y.$  Thus, $\as Y\le \as X.$
\end{proof}

\begin{Proposition} \label{R^n} $\as \mathbb{R}^n=n.$
\end{Proposition}

\begin{proof} We will see in the next section that $\as (X\times Y)\le
\as X+\as Y,$ or one could convince oneself of the upper bound by
drawing pictures for the plane and imagining their extensions to
higher dimensions.

For the lower bound, we use the technique of \cite{DKU}.  One could
also use homological methods, see \cite{Ro03}.  First, we assume
$\dim[0,1]^n=n$ to be known.  If $\as\mathbb{R}^n\le k,$ we can take
$k+1$ $R$-disjoint families of uniformly bounded subsets of
$\mathbb{R}^n$ which cover $\mathbb{R}^n.$ But, contracting the
covers and taking closures gives an $\epsilon$-cover of $[0,1]^n$
for small $\epsilon$ with multiplicity $\le k+1.$ Thus, $k\ge n.$
\end{proof}
\subsection{A union theorem for asdim}

In this section we establish a union theorem for asdim.  It should
be noted that here asymptotic dimension varies slightly from
covering dimension.  For example, the finite union theorem for
covering dimension says $\dim(X\cup Y)\le \dim X+\dim Y+1,$ and that
inequality is sharp.  Also, a countable union theorem for covering
dimension is: $\dim(\cup_i C_i)\le \max_{i}\{\dim C_i\}$ where the
$C_i$ are closed subsets of $X.$  Notice that there can be no direct
analog of this theorem for asymptotic dimension since every finitely
generated group is a countable set of points, and as we shall see,
finitely generated groups can have arbitrary (even infinite)
asymptotic dimension.

Let $\sU$ and $\sV$ be families of subsets of $X$.  Define the
$r$-saturated union of $\sV$ with $\sU$ by
\[\sV\cup_r\sU=\{N_r(V;\sU)\mid V\in\sV\}\cup \{U\in\sU\mid d(U,\sV)>r\},\]
where $N_r(V;\sU)=V\cup\bigcup_{d(U,V)\le r} U.$

\begin{Proposition} Let $\sU$ be an $r$-disjoint, $R$-bounded family
of subsets of $X$ with $R\ge r.$  Let $\sV$ be a $5R$-disjoint, $D$
bounded family of subsets of $X.$  Then $\sV\cup_r\sU$ is
$r$-disjoint and $(D+2(r+R))$-bounded.
\end{Proposition}

\begin{proof} The uniform bound on the diameters of elements of
$\sV\cup_r\sU$ is clear.  To see that the disjointness condition
holds, we consider the types of elements in $\sV\cup_r\sU,$ those of
the form $U$ and those of the form $N_r(V;\sU).$  Obviously if
$U\neq U'$ then $U$ and $U'$ are $r$-disjoint by the definition of
$\sU.$  A set of the form $U$ and a set $N_r(V;\sU)$ are
$r$-disjoint by definition of $N_r(V;\sU).$  Finally we consider
$N_r(V;\sU)$ and $N_r(V';\sU),$ where $V\neq V'.$  Clearly these
sets are contained in $N_{r+R}(V)$ and $N_{r+R}(V'),$ respectively.
Since $d(V,V')\ge 5R,$ and $R\ge r,$ we find that
$d(N_r(V;\sU),N_r(V';\sU))\ge r.$
\end{proof}

Let $X$ be a metric space.  We will say that the family
$\{X_\alpha\}$ of subsets of $X$ satisfies the inequality $\as
X_\alpha\le n$ {\em uniformly} if for every $r<\infty$ one can find
a constant $R$ so that for every $\alpha$ there exist $r$-disjoint
families $\sU^0_\alpha,\ldots,\sU^n_\alpha$ of $R$-bounded subsets
of $X_\alpha$ covering $X_\alpha.$ A typical example of such a
family is a family of isometric subsets of a metric space. Another
example is any family containing finitely many sets.

\begin{Theorem} Let $X=\cup_\alpha X_\alpha$ be a metric space where
the family $\{X_\alpha\}$ satisfies the inequality $\as X_\alpha\le
n$ uniformly.  Suppose further that for every $r$ there is a
$Y_r\subset X$ with $\as Y_r\le n$ so that
$d(X_\alpha-Y_r,X_{\alpha'}-Y_r)\ge r$ whenever $X_\alpha\neq
X_{\alpha'}.$  Then $\as X\le n.$
\end{Theorem}

Before proving this theorem, we state a corollary: the finite union
theorem for asymptotic dimension.

\begin{Corollary} Let $X$ be a metric space with $A,B\subset X$.
Then $\as (A\cup B)\le\max\{\as A, \as B\}.$
\end{Corollary}

\begin{proof}[Proof of Corollary.] Apply the union theorem to the
family $A,B$ with $B=Y_r$ for every $r.$
\end{proof}

\begin{proof}[Proof of Union Theorem.] Let $r<\infty$ be given and
take $r$-disjoint, $R$-bounded families $\sU^i_\alpha$ $(i=0,\ldots,
n)$ of subsets of $X_\alpha$ so that $\cup_i\sU^i_\alpha$ covers
$X_\alpha.$  We may assume $R\ge r.$ Take $Y=Y_{5R}$ as in the
statement of the theorem and cover $Y$ by families
$\sV^0,\ldots,\sV^n$ which are $D$-bounded and $5R$-disjoint. Let
$\bar\sU^i_\alpha$ denote the restriction of $\sU^i_\alpha$ to the
set $X_\alpha-Y.$  For each $i,$ take
$\sW^i_\alpha=\sV^i\cup_r\bar\sU^i_\alpha.$ By the proposition
$\sW^i$ consists of uniformly bounded sets and is $r$-disjoint.
Finally, put $\sW^i=\{W\in \sW^i_\alpha\mid \alpha\}$. Observe that
each $\sW^i$ is $r$-disjoint and uniformly bounded. Also, it is easy
to check that $\cup_i\sW^i$ covers $X.$
\end{proof}

\subsection{Connection to the classical dimension}

Let $\varphi\colon X\to\mathbb{R}$ be a function defined on a metric
space $X$. For every $x\in X$ and every $r>0$ let
$V_r(x)=\sup\{|\varphi(y)-\varphi(x)|\mid y\in N_r(x)\}$. A function
$\varphi$ is called {\it slowly oscillating\/} whenever for every
$r>0$ we have $V_r(x)\to0$ as $x\to\infty$ (the latter means that
for every $\varepsilon>0$ there exists a compact subspace $K\subset
X$ such that $|V_r(x)|<\varepsilon$ for all $x\in X\setminus K$).
Let $\bar X$ be the compactification of $X$ that corresponds to the
family of all continuous bounded slowly oscillating functions. The
{\it Higson corona\/} of $X$ is the remainder $\nu X=\bar X\setminus
X$ of this compactification.

It is known that the Higson corona is a functor from the category of proper
metric space and coarse maps into the category of compact Hausdorff spaces. In
particular, if $X\subset Y$, then $\nu X\subset \nu Y$.

For any subset $A$ of $X$ we denote by $A'$ its trace on $\nu X$,
i.e. the intersection of the closure of $A$ in $\bar X$ with $\nu
X$. Obviously, the set $A'$ coincides with the Higson corona $\nu
A$.

Dranishnikov, Keesling and Uspenskij \cite{DKU} proved the
inequality
$$
\dim\nu X\le\as X.
$$
It was shown there that $\dim\nu X\ge\as X$ for a large class of
spaces, in particular for $X=\mathbb{R}^n$. This gives another
approach to Proposition \ref{R^n}. Later Dranishnikov proved
\cite{Dr1} that the equality $\dim\nu X=\as X$ holds provided $\as
X<\infty$. The question of whether there is a metric space $X$ with
$\as X=\infty$ and $\dim\nu X<\infty$ is still open.

\subsection{Asymptotic inductive dimension}

The notion of asymptotic inductive dimension $\Ind$ was introduced
in \cite{Dr2}.

Recall that a closed subset $C$ of a topological space $X$ is a {\it
separator\/}
between disjoint subsets $A, B\subset X$ if $X\setminus C=U\cup V$, where $U,V$
are open subsets in $X$, $U\cap V=\emptyset$, $A\subset U$, $V\subset B$. A
closed subset $C$ of a topological space $X$ is a {\it cut\/}
between disjoint subsets $A, B\subset X$ if every continuum (compact connected
space) $T\subset X$ that intersect both $A$ and $B$ also intersects $C$.

 Let $X$ be a proper metric space. A subset
$W\subset X$ is called an {\it asymptotic neighborhood\/} of a subset $A\subset
X$ if $\lim_{r\to\infty}d(X\setminus N_r(x_0),X\setminus W)=\infty$. Two sets
$A,B$
in a metric space are {\it asymptotically disjoint\/} if $\lim
_{r\to\infty}d(A\setminus N_r(x_0),B\setminus N_r(x_0))=\infty$. In other words,
two
sets are asymptotically disjoint if the traces $A'$, $B'$ on $\nu X$  are
disjoint.

A subset $C$ of a metric space $X$ is an  {\it asymptotic separator\/}  between
asymptotically disjoint subsets $A_1,A_2\subset X$ if the trace $C'$ is a
separator
in $\nu X$ between $A_1'$ and $A_2'$.

We recall the definition of the asymptotic Dimensiongrad in the
sense of Brouwer $\Ind_b$ from \cite{DrZ}.

Let $X$ be a metric space and $\lambda>0$. A finite sequence $x_1,\dots,x_k$ in
$X$ is a {\it $\lambda$-chain\/} between subsets  $A_1,A_2\subset X$ if $x_1\in
A_1$, $x_k\in A_2$ and $d(x_i,x_{i+1})<\lambda$ for every $i=1,\dots, k-1$.
We say that a subset $C$ of a metric space $X$ is an  {\it asymptotic
cut\/}  between
asymptotically disjoint subsets $A_1,A_2\subset X$ if for every $D>0$ there is
$\lambda>0$ such that every $\lambda$-chain between $A_1$ and $A_2$ intersects
$N_D(C)$.

By definition, $\Ind X=\Ind_b X=-1$ if and only if $X$ is bounded.
Suppose we have defined the class of all proper metric spaces $Y$
with $\Ind Y\le n-1$ (respectively with $\Ind_b Y\le n-1$). Then
$\Ind X\le n$ (respectively  $\Ind_b Y\le n$) if and only if for
every asymptotically disjoint subsets $A_1,A_2\subset X$ there
exists an  asymptotic separator (respectively asymptotic cut) $C$
between $A_1$ and $A_2$ with $\Ind C\le n-1$ (respectively $\Ind_b
C\le n-1$). The dimension functions $\Ind$ and $\Ind_b$ are called
the {\em asymptotic inductive dimension\/} and {\em asymptotic
Brouwer inductive dimension\/} respectively.

It's easy to prove that $\Ind_bX\le\Ind X$, for
every $X$.
It is unknown if $\Ind=\Ind_b$ for proper metric spaces.

\begin{Theorem} For all proper metric spaces $X$ with $0<\as X<\infty$ we have
$$
\as X=\Ind X.
$$
\end{Theorem}

This theorem is a very important step in the existing proof of the
exact formula of the asymptotic dimension of the free product $\as
A\ast B$ of groups \cite{BDK}, see section 3.4.

Notice that there is a small problem with coincidence of $\as$ and
$\Ind$ in dimension $0.$  This leads to philosophical discussions of
whether bounded metric spaces should be defined to have $\as=-1$ or
$0.$ Observe that in the world of finitely generated groups, $\as
\Gamma=0$ if and only if $\Gamma$ is finite. On the other hand,
there are metric spaces, for instance $2^n\subset\mathbb{R}$ that
are unbounded yet have asymptotic dimension $0.$

\subsection{Embeddings into trees}

In \cite{Dr3}, Dranishnikov showed that every proper metric space
$X$ with $\as X\le n$ admits a uniform embedding into a product of
$n+1$ locally finite regular $\mathbb{R}$-trees. Using this result
Dranishnikov and Zarichnyi \cite{DrZ} constructed a metric space
$M_n$ with $\as M_n=n$ that is universal for the class of proper
metric spaces with $\as X\le n$. The space $M_n$ plays a crucial
role in the proof of the above theorem.

We recall the notion of Assouad-Nagata dimension
N-$\dim$ defined in  \cite{A} as follows: for a
metric space $X$ we have N-$\dim X\leq n$ if there
is a constant $C$ such that, for each $r>0$, there is an open cover $\sU(r)$ of
$X$ by sets of diameter $\le Cr$ such that
each open ball of radius $r$ meets at most $n+1$ members of $\sU(r)$.

U. Lang and T. Schlichenmaier gave the following refinement of
Dranishnikov's embedding \cite{La}:

\begin{Theorem} If for a metric space N-$\dim (X,d)\le n$ then
for sufficiently small $\epsilon$, $(X,d^{\epsilon})$ admits a
bi-Lipschitz embedding in the product of $n+1$ locally finite trees.
\end{Theorem}

Dranishnikov and Schroeder \cite{DS} proved that the hyperbolic
plane admits a bi-Lipschitz embedding into the product of two binary
trees.

Later Buyalo and Schroeder \cite{BuSch}, using techniques of
\cite{DS} and some results of Buyalo \cite{Bu} showed that a
finitely generated hyperbolic group $\Gamma$ can be
quasi-isometrically embedded into a product of $n$ binary trees
where $\dim\partial \Gamma=n$ and that this result is optimal.

Combining this result with work of \'{S}wi\k{a}tkowski \cite{Swi}
saying that $\as\Gamma\ge\dim\partial\Gamma+1$ for finitely
generated hyperbolic groups, one sees that
$\as\Gamma=\dim\partial\Gamma+1.$  Notice too that this is not the
case for hyperbolic spaces.  The Comb Space, $C$, of Papasoglu and
Gentimis \cite{PG} has $\dim \partial C=0$ whereas $\as C=2.$  This
example can be easily modified to give arbitrarily high asdim
without changing $\dim\partial C.$

\section{Hurewicz-type theorem for asdim}

In this section we state (without proof) a {\em Hurewicz Theorem}
for asymptotic dimension.  This theorem allows us to compute
asymptotic dimension in a myriad of situations including direct
products, free products of groups, and group extensions.

\subsection{Groups acting on metric spaces}

Before stating the asymptotic version of the Hurewicz theorem, we
state with proof an easier special case of the theorem for groups
acting on finite dimensional metric spaces.  This gives the flavor
of the proof of the Hurewicz theorem. In order to state the result,
we need the idea of an $R$-stabilizer.  It is a metric subspace of
the group $\Gamma$, not a subgroup.

Let $\Gamma$ act on the metric space $X$ by isometries. Let $R>0$ be
given. Let $x_0\in X$. Define the $R$-stabilizer of $x_0$ by
$W_R(x_0)=\{\gamma\in\Gamma\mid d(\gamma.x_0,x_0)\le R\}.$

\begin{Theorem} Assume that the finitely generated group $\Gamma$
acts by isometries on the metric space $X$ with $x_0\in X$ and $\as
X\le k.$  Suppose further that $\as W_R(x_0)\le n$ for all $R$. Then
$\as \Gamma<\infty.$
\end{Theorem}

The coarseness of the upper bound is a consequence of the use of
covers. As mentioned following the definition of asdim, proving a
tight upper bound requires the use of maps to uniform polyhedra.

\begin{proof} We will show that $\as \Gamma\le (n+1)(k+1)-1.$

We define a map $\pi:\Gamma\to X$ by the formula $\pi(g)=g(x_0)$.
Then $W_R(x_0)=\pi^{-1}(B_R(x_0))$. Let
$\lambda=\max\{d_X(s(x_0),x_0)\mid s\in S\}$. We show now that $\pi$
is $\lambda$-Lipschitz. Since the metric $d_S$ on $\Gamma$ is
induced from the geodesic metric on the Cayley graph, it suffices to
check that $d_X(\pi(g),\pi(g'))\le\lambda$ for all $g,g'\in\Gamma$
with $d_S(g,g')=1$. Without loss of generality we assume that
$g'=gs$ where $s\in S$. Then
$d_X(\pi(g),\pi(g'))=d_X(g(x_0),gs(x_0))=d_X(x_0,s(x_0))\le\lambda$.

Note that $\gamma B_R(x)=B_R(\gamma(x))$ and
$\gamma(\pi^{-1}(B_R(x)))= \pi^{-1}(B_R(\gamma(x)))$ for all
$\gamma\in\Gamma$, $x\in X$ and all $R$.

Given $r>0$, there are $\lambda r$-disjoint,  $R$-bounded families
${\mathcal F}^0,\dots,{\mathcal F}^k$ on the orbit $\Gamma x_0$. Let
 ${\sV}^0,\dots,{\sV}^n$ on $W_{2R}(x_0)$
 be $r$-disjoint uniformly bounded families given by the definition of
the inequality $\as W_R(x_0)\le n$. For every element $F\in{\mathcal
F}^i$ we choose an element $g_F\in\Gamma$ such that $g_F(x_0)\in F$.
We define $(k+1)(n+1)$ families of subsets of $\Gamma$ as follows
$$
{\mathcal W}^{ij}=\{g_F(C)\cap\pi^{-1}(F)\mid F\in{\mathcal F}^i,
C\in{\sV}^j\}
$$
Since  multiplication by $g_F$ from the left is an isometry, every
two distinct sets $g_F(C)$ and $g_F(C')$ are $r$-disjoint. Note that
$\pi(g_F(C)\cap\pi^{-1}(F))$ and $\pi(g_{F'}(C)\cap\pi^{-1}(F'))$
are $\lambda r$-disjoint for $F\neq F'$. Since $\pi$ is
$\lambda$-Lipschitz, the sets $g_F(C)\cap\pi^{-1}(F)$ and
$g_{F'}(C')\cap\pi^{-1}(F')$ are $r$-disjoint. The families
${\mathcal W}^{ij}$ are uniformly bounded, since the families
${\mathcal V}^j$ are, and  multiplication by $g$ from the left is an
isometry on $\Gamma$. We check that the union of the families
${\mathcal W}^{ij}$ forms a cover of $\Gamma$. Let $g\in\Gamma$ and
let $\pi(g)=F$, i.e. $g(x_0)\in F$. Since $\diam F\le R$, $x_0\in
g_F^{-1}(F)\le R$ and $g_F^{-1}$ acts as an isometry, we have
$g^{-1}_F(F)\subset B_R(x_0)$. Therefore, $g_F^{-1}g(x_0)\in
B_R(x_0)$, i. e. $g^{-1}_Fg\in W_R(x_0)$. Hence $g_F^{-1}g$ lies in
some set $C\in{\sV}^j$ for some $j$. Therefore $g\in g_F(C)$. Thus,
$g\in g_F(C)\cap\pi^{-1}(F)$.
\end{proof}

\subsection{Hurewicz-type theorem and applications}

In general it is not possible to say anything about the asymptotic
dimension of a surjective image based on the asymptotic dimension of
the domain space. All one has to do is check that a countable set
can be metrized as $\mathbb{Z}^n$ for any $n$ and so, as we shall
see, can be made to have arbitrary (even infinite) asymptotic
dimension.  On the other hand, when one has a map from a space $X$
and one knows the asymptotic dimension of the codomain, one can
often estimate $\as X.$  This is the situation of the following
theorem.

\begin{Theorem}[Hurewicz-type theorem] Let $f:X\to Y$ be a
Lipschitz map from a geodesic metric space to a metric space.
Suppose that for every $R<\infty$ the set $\as f^{-1}(B_R(y))\le n$
uniformly (in $y\in Y$). Then $\as X\le n+\as Y.$
\end{Theorem}

The proof of this theorem uses the definition of asymptotic
dimension in terms of Lipschitz maps to uniform complexes.  A
complete proof can be found in \cite{BD3}.  Another simpler version
of this result is the main theorem in \cite{BD2} involving the case
where $Y$ is a tree.

The Hurewicz theorem allows us to estimate the asymptotic dimension
of a direct product of metric spaces:

\begin{Corollary} Let $X$ and $Y$ be metric spaces.  Then $\as X\times
Y\le \as X+\as Y.$
\end{Corollary}

Although it is not difficult to prove that the product of two spaces
with finite asymptotic dimension has finite asymptotic dimension,
getting a sharp upper bound from definitions involving covers is
difficult, compare trying to pass from a cover of $\mathbb{R}$ to
one for $\mathbb{R}^2.$ However, with a little work, one can prove
this sharp upper bound for the $\as$ of a product from the
definition of $\as$ involving uniformly co-bounded Lipschitz maps to
uniform polyhedra.  Instead, we apply the Hurewicz-type theorem.

\begin{proof} The map $f:X\times Y\to Y$ given by $f(x,y)=y$ is
clearly Lipschitz since $d(x\times y, x'\times
y')=\left(d_X(x,x')^2+d_Y(y,y')^2\right)^{\frac12}\le d_X(x,x').$ It
remains only to show that $\as f^{-1}(B_R(y))\le \as X$ uniformly
for some $n.$ To see this, for any $r<\infty$ take $n+1$ families of
$r$-disjoint, $B$-bounded sets $\sU^0,\ldots,\sU^n$ whose union
covers $X$ and define $\sV^i_y=\{U\times B_R(y)\mid U\in \sV^i\}.$
For each $y\in Y,$ the $\sV^i_y$ are $r$-disjoint and
$\sqrt2\max\{B,R\}$-bounded.  Obviously, the union of the $\sV_y^i$
forms a cover of $f^{-1}(B_R(y))$ for each $y.$  Thus, the family
satisfies the inequality $\as f^{-1}(B_R(y))\le\as X$ uniformly.  By
the Hurewicz theorem, $\as(X\times Y)\le \as X+\as Y.$
\end{proof}

As another corollary of the Hurewicz-type theorem, we arrive at the
case of interest to geometric group theorists: a finitely generated
group acting by isometries on a metric space.

\begin{Corollary} \label{ActionOnSpace} Let $\Gamma$ be a finitely generated group acting
by isometries on a metric space $X$.  Fix some $x_0$ in $X$ and
suppose that $\as W_R(x_0)\le k$ for all $R$.  Then, $\as\Gamma\le
k+\as X.$
\end{Corollary}

\begin{proof} Fix a symmetric generating set $S$ for $\Gamma.$  Let
$\lambda=\max\{d_X(s.x_0,x_0)\mid s\in S\}.$  Define $\pi:\Gamma\to
X$ by $\pi(\gamma)=\gamma.x_0.$  We claim that $\pi$ is
$\lambda$-Lipschitz and that $\as\pi^{-1}(B_R(x))\le k$ uniformly in
$x\in \Gamma.x_0.$

Since $\Gamma$ is a finitely generated group with the word metric it
is a discrete geodesic space, so it suffices to check the Lipschitz
condition on pairs at distance $1$ from each other.  It is easy to
see that such a pair must be of the form $(\gamma, \gamma s),$ where
$s\in S.$  We compute: \[d_X(\pi(\gamma),\pi(\gamma
s))=d_X(\gamma.x_0,\gamma s.x_0)=d_X(x_0,s.x_0)\le \lambda,\] so
$\pi$ is $\lambda$-Lipschitz.  Finally, it is easy to check that
$\pi^{-1}(B_R(g.x_0))=gW_R(x_0).$  Since left multiplication is an
isometry, the sets $\pi^{-1}(B_R(g.x_0))$ are all isometric to
$W_R(x_0),$ so $\as\pi^{-1}(B_R(g.x_0))\le k$ uniformly.  Finally,
since $\as\Gamma.x_0\le \as X,$ we get the desired inequality.
\end{proof}

Another application of the Hurewicz theorem is to group extensions,
i.e. groups $G$ arising in exact sequences of the form: \[1\To K\To
G\To H\To 1.\]

\begin{Corollary} Let $f:G\to H$ be a surjective homomorphism of
finitely generated groups with $\ker f=K.$  Then $\as G\le \as H+\as
K.$
\end{Corollary}

\begin{proof} Let $S$ be a symmetric generating set for $G$ and take
$\bar S=f(S)$ to be a generating set for $H.$  If
$x^{-1}y=s_{i_1}\cdots s_{i_k}$ is a shortest presentation in terms
of generators, then $d(x,y)=k$ and
$d_H(f(x),f(y))=\|f(x^{-1}y)\|=\|\bar{s}_{i_1}\cdots\bar{s}_{i_k}\|\le
k.$  Thus, $f$ is $1$-Lipschitz.

Next we claim that $W_R(e)=N_R(K).$  Since $N_R(K)$ is
quasi-isometric to $K,$ this will say that $\as W_R(e)=\as K$ for
all $R$ so that we may apply the Hurewicz theorem to get the desired
result.

First suppose that $g\in W_R(e).$  Then $\|f(g)\|\le R.$ Thus,
$f(g)=\bar{s}_{i_1}\cdots\bar{s}_{i_k}$ with $k\le R.$  If
$u=s_{i_1}\cdots s_{i_k},$ then $gu^{-1}\in K$ and
$d_G(g,gu^{-1})=k\le R.$  On the other hand, if $d_G(x,K)\le R,$
then since $f$ is $1$-Lipschitz, $d_H(f(x),e)\le R.$
\end{proof}

\begin{Corollary} Let $\Gamma$ be a finitely generated polycyclic
group with Hirsch length $h(\Gamma)=n.$  Then $\as \Gamma\le n.$
\end{Corollary}

\begin{proof} Since $\Gamma$ is polycyclic, there exists a chain
\[1=\Gamma_0\lhd \Gamma_1\lhd \cdots\lhd\Gamma_n=\Gamma\]
where each $\Gamma_{i+1}/\Gamma_i$ is cyclic.  The Hirsch length is
$h(\Gamma)=\sum\text{rk}(\Gamma_{i+1}/\Gamma_i).$  Applying the
extension theorem we see that $\as\Gamma\le h(\Gamma).$
\end{proof}

Since every finitely generated nilpotent group is polycyclic we
immediately obtain the following result.

\begin{Corollary} \label{6} Let $\Gamma$ be a finitely generated nilpotent
group.  Then $\as\Gamma\le h(\Gamma).$
\end{Corollary}

\begin{example} Let $H$ denote the $3\times 3$ integral Heisenberg group.  Then
$\as H\le 3.$
\end{example}

Corollary \ref{6} can be extended to nilpotent Lie groups $N$ if one
defines the Hirsch length $h(N)$ as the sum of the number of factors
in $\Gamma_{i+1}/\Gamma_i$ isomorphic to $\mathbb{R}$ for the
central series $\{\Gamma_i\}$ of $N$. We take an equivariant metric
on $N$ and on the quotients. Then the projection $\Gamma_{i+1}\to
\Gamma_{i+1}/\Gamma_i$ is 1-Lipschitz and $\Gamma_{i+1}/\Gamma_i$ is
coarsely isomorphic to $\mathbb{R}^{n_i}$. Then we have

\begin{Corollary}  \label{Lie} Let $N$ be a nilpotent Lie group endowed with an equivariant
metric. Then $\as N\le h(N).$
\end{Corollary}

Since $h(N)=\dim N$ for simply connected $N$, we obtain

\begin{Corollary} \cite[Theorem 3.5]{CG}  \label{cg} For
a simply connected nilpotent Lie group $N$ endowed with an
equivariant metric $\as N\le \dim N.$
\end{Corollary}

Actually in view of \cite[Corollary 1.F1]{Gr93} the inequalities in
Corollaries \ref{Lie} and \ref{cg} are equalities.

Corollary \ref{cg} is the main step in the proof of the following

\begin{Theorem}\cite{CG} For a connected Lie group $G$ and its maximal compact subgroup $K$
there is a formula $\as G/K=\dim G/K$ where $G/K$ is endowed with a
G-invariant metric.
\end{Theorem}

\subsection{Asymptotic dimension of hyperbolic space}

This last theorem in particular allows to show that the asymptotic
dimension of the hyperbolic space $\mathbb{H}^n$ is $n$.

\begin{Corollary} $\as \mathbb{H}^n=n$.
\end{Corollary}

\begin{proof}
Take $G=O(n,1)_+$ and $K=O(n)$.
\end{proof}
This computation can be generalized in spirit of \cite{Ro03}.

Let $(X,d)$ be a metric space. By $\sH(X)$ we denote the space of
balls in $X$ endowed with the following metric
$$
\rho(B_t(x),B_s(y))=2\ln\left(\frac{d(x,y)+\max\{t,s\}}{\sqrt{ts}}\right).
$$

We note that $\sH(\mathbb{R}^n)$ is coarsely equivalent to
$\mathbb{H}^{n+1}$, \cite[Example 2.60]{Ro03}.

We recall that a  metric space $X$ with $\as X\le n$ is said to
satisfy the {\em Higson property} \cite{DrZ} if there exists $C>0$
such that for every $D>0$ there exists a cover $\sU$ of $X$ with
$\mesh(\sU)<CD$ and such that $\sU=\sU^0\cup\dots\cup\sU^n$, where
$\sU^0,\dots,\sU^n$ are $D$-disjoint. In \cite{Ro03} $X$ satisfying
this condition are said to have asymptotic dimension $\le n$ {\em of
linear type}. Note that this condition is equivalent to the
asymptotic inequality $N-\dim X\le n$ for the Assouad-Nagata
dimension. It is shown in \cite{DrZ} that every metric space of
bounded geometry with $\as X\le n$ admits a coarsely equivalent
metric with the Higson property. Unfortunately the coarse type of
$\sH(X)$ depends on a metric on $X$ not only the coarse class of
metrics.

\begin{Theorem} Suppose that the metric space $(X,d)$
possesses the Higson property. Then $\as\sH(X)=\as X+1$.
\end{Theorem}
\begin{proof}
Consider the projection $\pi:\sH(X)\to\mathbb{R}$ defined by
$\pi(B_t(x))=\ln t$ and apply the Hurewicz-type Theorem to it (see
\cite[Corollary 9.21]{Ro03}).
\end{proof}

\subsection{Groups acting on trees and Bass-Serre theory}

For the basics of Bass-Serre theory the reader is referred to
\cite{Se} or (for generalizations of Bass-Serre theory) to
\cite{BH}.  Recall that a tree $T$ has $\as T\le 1.$

The Bass-Serre theory allows us to consider various constructions
with groups simultaneously.  In particular, given two finitely
generated groups $A$ and $B$ we would like to be able to estimate
the asymptotic dimension of their free product, $A\ast B$ an
amalgamated free product $A\ast_C B$ or an HNN-extension formed from
one of these groups.  The Bass-Serre theory tells us that these are
all examples of groups which act co-compactly by isometries on
trees.

\begin{Corollary} Let $\Gamma$ be a finitely generated group acting
co-compactly by isometries on a tree $T$.  Suppose that for all
vertices $v,$ $\as\Gamma_v\le n.$ Then $\as\Gamma\le n+1.$
\end{Corollary}

Applying Corollary \ref{ActionOnSpace} we see that the only thing we
need to show is that for all $R$, $\as W_R(x)\le n.$ This is not
obvious.  For the proof, the reader is referred to \cite[Lemma
3]{BD2}.  Instead, we offer a simpler case, that of the free product
$A\ast B.$

\begin{Corollary} Let $A$ and $B$ be finitely generated groups with
$\as A\le n$ and $\as B\le n.$  Then $\as A\ast B\le n+1.$
\end{Corollary}

\begin{proof} The free product $A\ast B$ is the fundamental group of
the graph of groups with two vertices, labeled $A$ and $B$ and one
edge, labeled $\{e\}.$ This group acts on a tree by isometries so
that the quotient consists of two vertices and one edge.  The
vertices of the tree consist of formal cosets of either $A$ or $B$
in the group.  The vertices $xA$ and $yB$ are connected by an edge
in the tree when there is a $z\in A\ast B$ such that $zA=xA$ and
$zB=xB.$

It is not difficult to show in this case that $W_R(eA)$ consists of
alternating products of the form $AB\cdots ABA$ of length $R+1$ when
$R$ is even and $AB\cdots BAB$ when $R$ is odd.

To see that the asymptotic dimension of such products does not
exceed $n$, one applies the union theorem for asymptotic dimension
and induction.

In the case $R=0$ we have $W_0(eA)=A,$ and $\as A\le n$ by
assumption.  For the case $R>0$ we assume $R$ is even (if it is odd,
the proof is essentially the same). Write $AB\cdots ABA$ as
$\cup_{x\in AB\cdots AB} xA.$  By assumption $\as AB\cdots AB\le n$,
and since $xA$ is isometric to $A$ we know that $\as xA\le n$
uniformly.

To apply the union theorem it remains only to find $Y_r$ so that
$\as Y_r\le n$ and so that the sets $xA-Y_r$ and $x'A-Y_r$ are
$r$-disjoint when $xA\neq xA'.$  To this end, set $Y_r=AB\cdots
ABB_r(e)$ where $B_r(e)$ denotes the ball of radius $R$ around $e$
taken in $A.$ This is coarsely isometric to $AB\cdots AB$, which, by
the inductive hypothesis, has asymptotic dimension not exceeding
$n.$  Finally, if $xa$ and $x'a'$ are in distinct $xA-Y_r$ and
$x'A-Y_r$ then $d(xa,x'a')=\|a^{-1}x^{-1}x'a'\|\ge \|a\|+\|a'\|\ge
r.$
\end{proof}

Using the asymptotic inductive dimension, with Keesling the authors
were able to give an exact formula for the asymptotic dimension of
such a free product in \cite{BDK}: $\as A\ast B=\max\{n,1\}.$  An
exact formula in the case of amalgamated free products still does
not exist.

Applying the Bass-Serre theory, one can prove this result for the
more general situation of graphs of groups, (see \cite{BD2}).  This
situation includes amalgamated free products and HNN extensions as
special cases.

\begin{Theorem} Let $\pi$ be the fundamental group of a finite graph
of groups where all vertex groups satisfy $\as \Gamma_v\le n.$  Then
$\as \pi\le n+1.$
\end{Theorem}

Also, one can extend this result to complexes of groups in a natural
way, see \cite{Be05}.

\begin{Corollary} Let $\Gamma$ be a finitely generated group with one defining relator.
Then $\as\Gamma<\infty.$
\end{Corollary}

\begin{proof}An result of Moldavanskii from  1967 states that a finitely
generated one-relator group is an HNN extension of a finitely
presented group with shorter defining relator or is cyclic, see
\cite{LS}. If $A$ is a finitely generated group with $\as A=n$ and
$A\ast_C$ is an HNN extension of $A,$ then since $A\ast_C$ is the
fundamental group of a loop of groups, the Hurewicz-type theorem
tells us that $\as A\ast_C\le n+1.$  Iterating this procedure
finitely many times gives the desired result.
\end{proof}

The class of limit groups consists of those groups which naturally
arise in the study of solutions to equations in finitely generated
groups.  One definitions is the following:  A finitely presented
group $L$ is a {\em limit group} if for each finite subset
$L_0\subset L$ there is a homomorphism to a free group which is
injective on $L_0.$  For more information the reader is referred to
\cite{BF} and the references therein.

The next result does not (to our knowledge) appear in the
literature.  It was pointed out to the first author by Bestvina to
be an easy consequence of the Hurewicz-type theorem and deep results
in the theory of limit groups.

\begin{Proposition} Let $L$ be a limit group.  Then $\as L<\infty.$
\end{Proposition}

\begin{proof}  There are two ways to see this.  The first way is to construct
$L$ (say with height h) via fundamental groups of graphs of groups
where vertices have height ($h-1$) and height 0 groups are free
groups, free abelian groups and surface groups.

Another is to combine Osin's theorem \cite{Os} that groups that are
hyperbolic relative to a collection of subgroups with finite
asymptotic dimension themselves have finite asymptotic dimension
with a result of Dahmani \cite{Dah} stating that limit groups are
relatively hyperbolic with respect to their maximal abelian
non-cyclic subgroups.
\end{proof}

\section{Motivation}

Although $\as$ was introduced in 1993, it did not garner much
attention until a paper of G. Yu in 1998 \cite{Yu1} in which he
proved that the Novikov higher signature conjecture holds for groups
with finite asymptotic dimension.  Below we list some related
results.

\begin{Theorem} \cite{Yu1} Let $\Gamma$ be a finitely generated
group with finite $B\Gamma$ and $\as\Gamma<\infty.$ Then the Novikov
Conjecture holds for $\Gamma.$
\end{Theorem}

Later, G. Yu \cite{Yu2} would generalize this theorem to the
following:

\begin{Theorem} \cite{Yu2} Let $\Gamma$ be a discrete metric space
with bounded geometry admitting a uniform embedding into Hilbert
space.  Then the coarse Baum-Connes conjecture holds for $\Gamma.$
\end{Theorem}

An easy way to verify that a discrete metric space admits a uniform
embedding into Hilbert space is to verify that it has {\em Property
A} \cite{Yu2}.  A discrete metric space $Z$ has property A if there
exist maps $\{a_n\}_{n\in\mathbb{N}}$, $a^n:Z\to \Prob(Z)$ such that
the following two conditions hold:
\begin{enumerate}
    \item for every $n$ there is an $R$ so that for every $z\in Z,$
    $\supp(a^n_z)\subset B_R(z)$; and
    \item for every $K>0$ \[\lim_{n\to\infty}\sup_{d(z,w)<K}\|a^n_z-a^n_w\|_1=0.\]
\end{enumerate}  Here, by way of notation, $a^n_z(x)=(a^n(z))(x).$

Before moving on to groups with finite asymptotic dimension, we give
an example.

\begin{example}  Let $T$ be a tree, then $T$ has property $A.$
\end{example}

\begin{proof} Fix some geodesic ray $\gamma$ in $T.$  It is easy to
see that for any point, $x\in T$ there is exactly one geodesic ray
from $x$ whose intersection with $\gamma$ is also a geodesic ray.
Let $\gamma_x$ be this unique ray for $x$.  We abuse notation by
writing $\gamma_x:\mathbb{N}\to T$ for the function yielding
$\gamma_x.$  For each $n,$ define $a^n_x=\frac{1}{n+1}\sum_{i=0}^n
\gamma_x(i)\delta_{\gamma_x(i)}.$

Clearly the support of $a^n_x$ is contained in the $n$ neighborhood
of $x$ and it is easy to check that for any $K$, $\|a^n_z-a^n_w\|\le
\frac{n-(n-k)}{n+1}$ which goes to $0$ as $n\to\infty.$
\end{proof}

To see that finitely generated groups with finite asymptotic
dimension admit a uniform embedding into Hilbert space we check that
they have property A. This result and its proof are due to Higson
and Roe, \cite{HR}, see also Theorem \ref{PropA}.

\begin{Theorem} \label{HR} Let $\Gamma$ be a finitely generated group
with finite asymptotic dimension.  Then $\Gamma$ has property A.
\end{Theorem}

\begin{proof} Suppose $\as \Gamma=n$ and $R<\infty$ is given.  Let
$\sU$ be a cover of $\Gamma$ by uniformly bounded sets with
multiplicity $\le n+1$ and Lebesgue number $>R.$  There is a
partition of unity $\{\phi\}$ subordinate to this cover which has
the following properties:
\begin{enumerate}
    \item each $\phi$ is Lipschitz with Lipschitz constant $<2/R$;
    \item $\sup_\phi(\diam(\supp\phi))<\infty$; and
    \item for any $\gamma$, at most $n+1$ of the $\phi(\gamma)\neq
    0.$
\end{enumerate}

Now, for any $\gamma\in \Gamma,$ take $\gamma_\phi$ to be a non-zero
value of $\phi.$  Then, put $a^R_\gamma=\sum_\phi
\phi(\gamma)\delta_{\gamma_\phi}.$

Each $a^R_\gamma$ is in $\Prob(\Gamma),$ condition (2) implies that
for each $R$ there is a $D$ so that $\supp(a^R_\gamma)\subset
B_D(\gamma)$ and conditions (1) and (3) say that for all $K$
\[\lim_{R\to\infty}\sup_{d(\gamma,\gamma')<K}\|a^R_\gamma-a^R_{\gamma'}\|_1=0\]
\end{proof}

Other implications of finite asymptotic dimension are the following:

\begin{Theorem} \cite{Dr2} If $\Gamma$ is
the fundamental group of an aspherical manifold $M$ and
$\as\Gamma<\infty$ then the universal cover of $M$ is hyperspherical
\end{Theorem}

\begin{Theorem} \cite{DFW} If $\Gamma$ is a finitely generated group
with finite $B\Gamma$ and $\as\Gamma<\infty$ then the integral Novikov conjecture holds for
$\Gamma.$
\end{Theorem}

\begin{Theorem} \cite{Ba,CG} If $\Gamma$ is a finitely generated group
with finite $B\Gamma$ and $\as\Gamma<\infty$ then the integral K-theoretic Novikov
conjecture holds for $\Gamma.$
\end{Theorem}

\section{Dimension growth}

Using the fact that $\as\mathbb{Z}^n=n$ it is not difficult to
construct an example of a finitely generated group with infinite
asymptotic dimension.  In fact, the reduced wreath product of
$\mathbb{Z}$ by $\mathbb{Z}$ (denoted $\mathbb{Z}\wr\mathbb{Z}$) is
finitely generated and contains a copy of $\mathbb{Z}^n$ for every
$n$, (see \cite{Ro03} for details). Therefore,
$\as\mathbb{Z}\wr\mathbb{Z}=\infty.$

Note that this result implies that even the finiteness of asymptotic
dimension is not preserved under quotients as this group is a
quotient of $\mathbb{F}_2$ and $\as\mathbb{F}_2=1.$

The invariant we associate to groups with infinite asymptotic
dimension is the {\em dimension growth} of the space.  We define the
dimension function of the metric space $X$ as follows:
\[d_X(\lambda)=\min\{m(\sU)-1\mid L(\sU)\ge\lambda,
\sup_{U\in\sU}\diam(U)<\infty, \sU\text{ covers }X\},\]
where $m(\sU)$ denotes the multiplicity of the cover $\sU$ and
$L(\sU)$ denotes a Lebesgue number for the cover.

Observe that $\lim_{\lambda\to\infty} d_X(\lambda)=\as X$ and that
the function is monotonic.

Obviously changing the metric on $X$ could drastically alter the
function $d_X,$ however, the growth of the function $d_X$ is an
invariant of quasi-isometry.

\begin{Proposition} \cite{Be7} Let $X$ and $Y$ be discrete metric spaces with
bounded geometry. Suppose that $X$ and $Y$ are quasi-isometric. Then
there is some positive $k$ so that $d_X(\lambda)\le
kd_Y(k\lambda+k)+k.$ In particular, the growth rate of $d_\Gamma$ is
well-defined for finitely generated group $\Gamma.$
\end{Proposition}

The proof is technical so it is omitted.

It is easy to see that the dimension function for a group cannot
grow arbitrarily fast:

\begin{Proposition} \cite{Dr04} For a finitely generated group
$\Gamma$ we have $d_\Gamma(\lambda)\le e^{\alpha\lambda}.$
\end{Proposition}

\begin{proof} There is some $a$ so that
$\vol(B_\lambda(x))<e^{a\lambda}.$  Take
$\sU_\lambda=\{B_\lambda(x)\mid x\in \Gamma\}.$  Then,
$L(\sU_\lambda)=\lambda$ and $m(U_\lambda)\le e^{a\lambda}.$
\end{proof}

Dranishnikov was able to generalize Theorem \ref{HR} of \cite{HR} to
spaces with bounded geometry where the growth of the dimension
function is linear \cite{Dr1}.  The best result along these lines is
the following:

\begin{Theorem} \cite{Dr04} \label{PropA} Let $\Gamma$ be a finitely generated
group with $d_\Gamma(\lambda)\le \lambda^m$.  Then $\Gamma$ has
property A; in particular, the Novikov higher signature conjecture
holds for $\Gamma.$
\end{Theorem}

We saw at the beginning of this section that
$\mathbb{Z}\wr\mathbb{Z}$ has infinite asdim.  The next result tells
us about the growth of its dimension function.  In particular this
gives an example of a group with infinite asdim whose dimension
function grows at most polynomially.

\begin{Proposition} Let $N$ be a finitely generated nilpotent group.
Suppose that $\as G<\infty$ where $G$ is a finitely generated group.
Then $d_{G\wr N}(\lambda)\le \lambda^n$ for some $n.$  Here, as
before, $G\wr N$ denotes the reduced wreath product.
\end{Proposition}

\providecommand{\bysame}{\leavevmode\hbox
to3em{\hrulefill}\thinspace}
\providecommand{\MR}{\relax\ifhmode\unskip\space\fi MR }
\providecommand{\MRhref}[2]{%
  \href{http://www.ams.org/mathscinet-getitem?mr=#1}{#2}
} \providecommand{\href}[2]{#2}

\end{document}